# CORRECTION

## SDES WITH OBLIQUE REFLECTIONS ON NONSMOOTH DOMAINS

*The Annals of Probability* **21** (1993) 554–580

By Paul Dupuis[1] and Hitoshi Ishii[2]

*Brown University and Waseda University*

It has been pointed out by Weining Kang and Ruth Williams that there is an error in an argument in [1]. The purpose of this note is to correct the argument.

The error affects only Case 2 of the paper, and occurs in the first display at the top of page 580, at the end of the proof of Theorem 5.1. This display claims that a certain bound follows from (3.28) of the paper, and implicitly assumes that if

$$(p, q) \doteq D f_\varepsilon^\beta (Y(s), Y'(s))$$

and if $Y(s) \in \partial G$ and $Y'(s) \in \overline{G}$, then $Y(s) + \beta p \in \partial G$ and $Y'(s) + \beta q \in \overline{G}$, which need not be true. The statement of the theorem is still correct, and the reason is that underlying assumptions are in some sense robust with respect to small perturbations of the boundary.

Before presenting the correction we review the assumption made in Case 2. There is an open set $W$ containing $\overline{G}$ and a $C^{2,+}$ function $g$ on $W \times \mathbb{R}^N$ which for each fixed $x \in W$ is $C^1$ as a mapping $r \mapsto g(x, r)$. Furthermore there are constants $C > 0$ and $\theta > 0$ such that

(3.21) $\qquad g(x, 0) = 0,$

(3.22) $\qquad g(x, r) \geq |r|^2,$

(3.23) $\qquad \langle D_r g(x, r), \gamma_i(x) \rangle \geq 0 \qquad \text{if } \langle r, n_i(x) \rangle \geq -\theta |r|,$

(3.24) $\qquad |p| \leq C |r|^2, \qquad |q| \leq C|r| \qquad \text{if } (p, q) \in D^+ g(x, r),$

Received October 2007; revised October 2007.
[1]Supported in part by NSF Grant DMS-04-04806.
[2]Supported in part by JSPS Grant-In-Aids for Scientific Research 15340051.

and for any $x \in \overline{G}, r \in \mathbb{R}^N$, there is $(p, q) \in D^+g(x, r)$ such that

$$(3.25) \qquad \left((p, q), C\begin{pmatrix} |r|^2 I & 0 \\ 0 & I \end{pmatrix}\right) \in D^{2,+}g(x, r).$$

A more careful statement of condition (3.23) is

$$\langle D_r g(x, r), \gamma_i(x) \rangle \geq 0 \qquad \text{if } x \in \partial G, \ i \in I(x) \text{ and } \langle r, n_i(x) \rangle \geq -\theta |r|.$$

The next lemma shows that (3.23) is in some sense robust.

LEMMA 1. *For each $R > 0$ there is a function $\sigma_R \in C([0, \infty))$, with $\sigma_R(0) = 0$, such that for any $x \in \partial G$, $i \in I(x)$, $y \in W$, and $r \in \mathbb{R}^N$, if $|r| \leq R$ and $\langle r, n_i(x) \rangle \geq -\theta |r|$, then*

$$\langle D_r g(y, r), \gamma_i(y) \rangle \geq -\sigma_R(|y - x|).$$

PROOF. It is enough to show the following for each $R < \infty$. For any $x_k \in \partial G$, $y_k \in W$, $i_k \in I(x_k)$ and $r_k \in \mathbb{R}^N$, if $|y_k - x_k| \to 0$, $|r_k| \leq R$ and $\langle r_k, n_{i_k}(x_k) \rangle \geq -\theta |r_k|$, then $\liminf_k \langle D_r g(y_k, r_k), \gamma_{i_k}(y_k) \rangle \geq 0$. Suppose this statement is false. Then there are $\eta > 0$ and $R > 0$ with the following property: for any $k \in \mathbb{N}$ there are $x_k \in \partial G$, $i_k \in I(x_k)$, $y_k \in W$, and $r_k \in \mathbb{R}^N$ such that $|y_k - x_k| < 1/k$, $|r_k| \leq R$, $\langle r_k, n_{i_k}(x_k) \rangle \geq -\theta |r_k|$, and

$$\langle D_r g(y_k, r_k), \gamma_{i_k}(y_k) \rangle < -\eta.$$

We may assume that as $k \to \infty$,

$$x_k \to x \in \partial G, \qquad r_k \to r.$$

Since $|x_k - y_k| < 1/k$, we have

$$y_k \to x \qquad \text{as } k \to \infty.$$

We may assume as well that $i_k = i$ for all $k$ and for some $i \in I$.

Set $q_k = D_r g(y_k, r_k)$ and for each $k$ choose $p_k \in \mathbb{R}^N$ so that $(p_k, q_k) \in D^+g(y_k, r_k)$. By (3.24) we have $|p_k| \leq CR^2$ and $|q_k| \leq CR$. Thus we may assume that, as $k \to \infty$, $(p_k, q_k) \to (p, q) \in \mathbb{R}^{2N}$. By the semiconcavity of $g$, we have $(p, q) \in D^+g(x, r)$. In particular, $q = D_r g(x, r)$. Since

$$\langle q_k, \gamma_i(y_k) \rangle < -\eta \qquad \text{for all } k,$$

we get

$$\langle q, \gamma_i(x) \rangle \leq -\eta.$$

Also, since

$$\langle r_k, n_i(x_k) \rangle \geq -\theta |r_k| \qquad \text{for all } k,$$

we have
$$\langle r, n_i(x)\rangle \geq -\theta |r|.$$
Finally, since $i \in I(x_k)$ for all $k$, by the upper semicontinuity of $I$ we have $i \in I(x)$.

Thus we have
$$x \in \partial G, \qquad i \in I(x), \qquad \langle n_i(x), r\rangle \geq -\theta|r|,$$
$$\langle \gamma_i(x), q\rangle \leq -\eta, \qquad q = D_r g(x,r),$$
which contradicts (3.23). □

We now introduce an approximation to $g$ based on sup-convolution. Define $\varphi(x,r) \doteq (|x|^2 + A)|r|^2$. Since $G$ is bounded, if $A > 0$ is large enough, then for some $0 < \delta < C < \infty$
$$\delta \begin{pmatrix} |r|^2 I & 0 \\ 0 & I \end{pmatrix} \leq D^2\varphi(x,r) \leq C \begin{pmatrix} |r|^2 I & 0 \\ 0 & I \end{pmatrix}$$
for all $x \in G$.

Let
$$h(x,r) \doteq g(x,r) - \psi(x,r),$$
where $\psi = B\varphi$ and $B > 0$ is large enough that $h$ is concave. For $\beta \in (0,1)$ define
$$h^\beta(x,r) \doteq \sup\left\{h(y,s) - \frac{1}{2\beta}(|x-y|^2 + |r-s|^2)\right\}.$$
As is well known (and easy to check), the concavity of $h$ implies that $h^\beta$ is also concave. Finally, set
$$g^\beta(x,r) \doteq h^\beta(x,r) + \psi(x,r).$$
Since $h^\beta$ is concave,
$$((p,q), BD^2\psi(x,r)) \in D^{2,+}g^\beta(x,r)$$
for $(p,q) = Dg^\beta(x,r)$, and therefore
$$\left((p,q), BC\begin{pmatrix} |r|^2 I & 0 \\ 0 & I \end{pmatrix}\right) \in D^{2,+}g^\beta(x,r).$$
We will use $C$ for a constant that takes values in $(0,\infty)$ and whose value may change from line to line, but in all cases $C$ can be chosen so that it is independent of $\beta$. For $\varepsilon \in (0,1)$, let $f_\varepsilon^\beta(x,y) \doteq \varepsilon g^\beta(x,(x-y)/\varepsilon)$. It follows that for some $C \in (0,\infty)$,
$$D^2 f_\varepsilon^\beta \leq \frac{C}{\varepsilon}\begin{pmatrix} I & -I \\ -I & I \end{pmatrix} + \frac{C|x-y|^2}{\varepsilon}\begin{pmatrix} I & 0 \\ 0 & I \end{pmatrix}.$$

This is a key property required of the mollification $f_\varepsilon^\beta$ in calculations prior to page 580. It remains to show how the argument on page 580 can be replaced. We will use that

$$|D_x\psi(x,r)| \leq C|r|^2, \qquad |D_{xx}\psi(x,r)| \leq C|r|^2,$$

$$|D_r\psi(x,r)| \leq C|r|, \qquad |D_{xr}\psi(x,r)| \leq C|r|,$$

$$|D_{rr}\psi(x,r)| \leq C.$$

LEMMA 2. *There is a constant $\tau > 0$ with the following property. For each $R > 0$ there is a function $\omega_R \in C([0,\infty))$ with $\omega_R(0) = 0$, such that*

$$\langle D_r g^\beta(x,r), \gamma_i(x) \rangle \geq -\omega_R(\beta)$$

*if $x \in \partial G$, $i \in I(x)$, $|r| \leq R$ and $\langle r, n_i(x) \rangle \geq -\tau|r|$.*

PROOF. Assume that

$$R > 0, \qquad x \in \partial G, \qquad i \in I(x), \qquad |r| \leq R, \qquad \langle r, n_i(x) \rangle \geq -\tau|r|,$$

where $\tau \in (0, 1 \wedge [\theta/3])$. Let $(p,q) = Dg^\beta(x,r)$. If

$$(\hat{p}, \hat{q}) \doteq (p - D_x\psi(x,r), q - D_r\psi(x,r)) = Dh^\beta(x,r),$$

then

$$(\hat{p}, \hat{q}) \in D^+ h(x + \beta\hat{p}, r + \beta\hat{q}) = D^+ g(\hat{x}, \hat{r}) - D\psi(\hat{x}, \hat{r}),$$

where $(\hat{x}, \hat{r}) \doteq (x + \beta\hat{p}, r + \beta\hat{q})$. Hence,

$$(\hat{p} + D_x\psi(\hat{x}, \hat{r}), \hat{q} + D_r\psi(\hat{x}, \hat{r})) \in D^+ g(\hat{x}, \hat{r}).$$

Since (3.24) states that

$$|\xi| \leq C|s|^2, \qquad |\eta| \leq C|s| \qquad \text{for } (\xi, \eta) \in D^+ g(y,s),$$

we have

$$|\hat{p} + D_x\psi(\hat{x}, \hat{r})| \leq C|\hat{r}|^2, \qquad |\hat{q} + D_r\psi(\hat{x}, \hat{r})| \leq C|\hat{r}|.$$

This implies

$$\begin{aligned}
|\hat{q}| &\leq |\hat{q} + D_r\psi(\hat{x}, \hat{r})| + |D_r\psi(\hat{x}, \hat{r})| \\
&\leq C|\hat{r}| + C|\hat{r}| \\
&= C|\hat{r}| \\
&= C|r + \beta\hat{q}| \\
&\leq C|r| + \beta C|\hat{q}|.
\end{aligned}$$

Choosing $\beta > 0$ small enough, we may assume that $\beta C \leq \min\{1/2, \tau\}$, so that

$$|q - D_r\psi(x,r)| = |\hat{q}| \leq C|r|, \qquad \beta|\hat{q}| \leq \beta C|\hat{r}| \leq \tau|\hat{r}|.$$

We then obtain

$$|r| = |\hat{r} - \beta\hat{q}| \leq (1+\tau)|\hat{r}|,$$
$$|\hat{r}| = |r + \beta\hat{q}| \leq (1+\tau)|r|,$$
$$|q| \leq |q - D_r\psi(x,r)| + |D_r\psi(x,r)| \leq C|r| + C|r| = C|r|.$$

Also, we have

$$|\hat{p}| \leq |\hat{p} + D_x\psi(\hat{x},\hat{r})| + |D_x\psi(\hat{x},\hat{r})| \leq C|\hat{r}|^2 + C|\hat{r}|^2 \leq C|r|^2.$$

For later use, note that $|p| = |\hat{p} + D_x\psi(x,r)| \leq C|r|^2$.

Now, we compute that

$$\langle n_i(x), \hat{r}\rangle = \langle n_i(x), r\rangle + \langle n_i(x), \hat{r} - r\rangle$$
$$\geq -\tau|r| - |\hat{r} - r|$$
$$\geq -\tau(|\hat{r}| + \beta|\hat{q}|) - \beta|\hat{q}|$$
$$\geq -\tau|\hat{r}| - 2\beta|\hat{q}|$$
$$\geq -3\tau|\hat{r}|.$$

Since $3\tau \leq \theta$, we have

$$\langle n_i(x), \hat{r}\rangle \geq -\theta|\hat{r}|.$$

Note that $|\hat{r}| \leq (1+\tau)|r| \leq 2|r| \leq 2R$ and $|\hat{x} - x| = \beta|\hat{p}| \leq \beta C|r|^2$. Thus by Lemma 1, we have

$$\langle \gamma_i(\hat{x}), \hat{q} + D_r\psi(\hat{x},\hat{r})\rangle \geq -\sigma_{2R}(\beta C R^2).$$

Finally, we compute that

$$\langle \gamma_i(x), q\rangle$$
$$\geq \langle \gamma_i(\hat{x}), q\rangle - C|x - \hat{x}||q|$$
$$\geq \langle \gamma_i(\hat{x}), \hat{q} + D_r\psi(\hat{x},\hat{r})\rangle - C|D_r\psi(\hat{x},\hat{r}) - D_r\psi(x,r)| - \beta C|\hat{p}||r|$$
$$\geq -\sigma_{2R}(\beta C R^2) - C|\hat{r}||\hat{x} - x| - C|\hat{r} - r| - C\beta|r|^3$$
$$\geq -\sigma_{2R}(\beta C R^2) - C(\beta|r|^3 + \beta|r|)$$
$$\geq -\sigma_{2R}(\beta C R^2) - \beta C(R^3 + R).$$

Thus Lemma 2 is valid with $\omega_R(t) = \sigma_{2R}(CR^2 t) + C(R + R^3)t$. □

6 P. DUPUIS AND H. ISHIIWe have $D_x f_\varepsilon^\beta(x,y) = \varepsilon D_x g^\beta(x, (x-y)/\varepsilon) + D_r g^\beta(x, (x-y)/\varepsilon)$, and by a calculation in Lemma 2 $|D_x g^\beta(x, (x-y)/\varepsilon)| \leq C|x-y|^2/\varepsilon^2$. In the proof of Theorem 5.1, we choose $R = \mathrm{diam}(G)/\varepsilon$, and then replace the second to fourth lines on page 580 by the following:

$$E \int_0^t u(Y, Y') \langle D_x f_\varepsilon^\beta(Y, Y'), \gamma(Y) \rangle \, d|k|(s)$$
$$\leq CE \int_0^t u(Y, Y') \frac{|Y-Y'|^2}{\varepsilon} \, d|k|(s) + C\omega_R(\beta) E|k|(T).$$

There is now no problem in taking the limit as $\beta \to 0$, and the proof proceeds as before.

## REFERENCE

[1] Dupuis, P. and Ishii, H. (1993). SDEs with oblique reflections on nonsmooth domains. *Ann. Probab.* **21** 554–580. MR1207237Lefschetz Center for Dynamical Systems
Division of Applied Mathematics
Brown University
Providence, Rhode Island 02912
USA
E-mail: dupuis@dam.brown.edu

Department of Mathematics
School of Education
Waseda University
Tokyo 169-8050
Japan